\documentstyle[12pt,amssymb,twoside]{amsart} 
\advance \topmargin by -\headheight 
\advance \topmargin by -\headsep 
\newtheorem{lemma}{Lemma}[section]
\newtheorem{prop}[lemma]{Proposition}
\newtheorem{theorem}[lemma]{Theorem}
\title{On certain isomorphisms between absolute Galois groups} 
\author{M.Rovinsky} 
\begin{document} 
\maketitle 
\section{Introduction}
Let $k$ be an algebraically closed field of characteristic zero, 
$L$ its finitely generated extension 
of transcendence degree $\ge 2$, and $L'$ another finitely generated 
extension of $k$. It is a result of Bogomolov \cite{bo-two} that 
any isomorphism between ${\rm Gal}(\overline{L}/L)$ and 
${\rm Gal}(\overline{L'}/L')$ is induced by an isomorphism of fields 
$\overline{L}\longrightarrow\overline{L'}$ identifying $L$ with $L'$. 

If the transcendence degree of $L$ over $k$ is one, the group 
${\rm Gal}(\overline{L}/L)$ is free, and therefore, its structure tells 
nothing about the field $L$. 

Let $F$ be an algebraically closed extension of $k$ of transcendence 
degree one, and $G=G_{F/k}$ be the group of automorphisms over $k$ 
of the field $F$. Let the set of subgroups $U_L:={\rm Aut}(F/L)$ for 
all subfields $L$ finitely generated over $k$ be the basis of 
neighborhoods of the unity in $G$. 

Let $\lambda$ be a continuous automorphism of $G$.
The purpose of this note is to show that if $\lambda$ induces 
an isomorphism ${\rm Gal}(F/L)\stackrel{\sim}{\longrightarrow}
{\rm Gal}(F/L')$ then the fields $L$ and $L'$ are isomorphic 
(see Theorem \ref{auto-gen} below for a more precise statement). 

\subsection{Notations} For a field $F_1$ and its subfield $F_2$ we 
denote by $G_{F_1/F_2}$ the group of automorphisms of the field $F_1$ 
over $F_2$. Throughout the note $k$ is an algebraically closed field 
of characteristic zero, $F$ its algebraically closed extension of 
transcendence degree $1\le n<\infty$ and $G=G_{F/k}$. 
If $K$ is a subfield of $F$ then $\overline{K}$ denotes 
its algebraic closure in $F$. 

For a topological group $H$ we denote by $H^{\circ}$ its subgroup 
generated by the compact subgroups, and by $H^{{\rm ab}}$ 
the quotient of $H$ by the closure of its commutant. 

For a smooth projective curve $C$ over a field, 
${\rm Pic}^{\ge m}(C)$ (resp., ${\rm Pic}^{>m}(C)$) 
is the submonoid in ${\rm Pic}(C)$ of 
sheaves of degree $\ge m$ (resp., $>m$). 

\section{A Galois-type correspondence} 
We consider a topology on $G$ with the basis of neighborhood of an 
automorphism $\sigma:F\stackrel{\sim}{\longrightarrow}F$ over $k$ 
given by the cosets of the form $\sigma U_L$ for all subfields $L$ 
of $F$ finitely generated over $k$, where $U_L={\rm Aut}(F/L)$. 
This topology was introduced in \cite{p-s-sha}. 

One checks that the topology is Hausdorff, locally 
compact, and totally decomposable. 
\begin{prop}[\cite{p-s-sha}, Lemma 1, Section 3] 
\label{main-psha}  The map $$\{\mbox{{\rm subfields in} $F$ 
{\rm over} $k$}\}\longrightarrow
\{\mbox{{\rm closed subgroups in} $G$}\}\quad\mbox{given by}
\quad K\longmapsto{\rm Aut}(F/K)$$ is injective and induces 
one-to-one correspondences 
\begin{itemize} \item $\{\mbox{{\rm subfields} $K$ {\rm of} $F$ 
{\rm with} $k\subseteq K$ {\rm and} $F=\overline{K}$}\}\leftrightarrow
\{\mbox{{\rm compact subgroups of} $G$}\}$; 
\item $\left\{\begin{array}{c} \mbox{{\rm subfields} $K$ {\rm of} 
$F$ {\rm finitely}} \\ \mbox{{\rm generated over} $k$ {\rm with} 
$F=\overline{K}$} \end{array} \right\}\leftrightarrow
\{\mbox{{\rm compact open subgroups of} $G$}\}$. 

The inverse correspondences are given by 
$G\supset H\longmapsto F^H\subseteq F$. \hfill $\Box$ 
\end{itemize} \end{prop}

Denote by $G^{\circ}$ the subgroup of $G$ generated by the 
compact subgroups. Obviously, $G^{\circ}$ is an open normal 
subgroup in $G$. 

\section{Decomposition subgroups in abelian quotients} 
Let $n=1$. We are going to show that for any 
continuous automorphism $\lambda$ of $G$ and any $L$ of finite 
type over $k$ one has $\lambda(U_L)=U_{L'}$ for some $L'$ isomorphic to $L$. 

To do that we first need to construct decomposition 
subgroups in the abelian quotients $U_L^{{\rm ab}}$. 

For a smooth projective model $C$ of $L$ over $k$ set 
$\Phi_L={\rm Hom}({\rm Div}^0(C),\widehat{{\Bbb Z}}(1))$. 
By Kummer theory, $U_L^{{\rm ab}}={\rm Hom}
(L^{\times},\widehat{{\Bbb Z}}(1))$, so, as ${\rm Pic}^0(C)$ 
is a divisible group, the short exact sequence 
$1\longrightarrow L^{\times}/k^{\times}\longrightarrow{\rm Div}^0(C)
\longrightarrow{\rm Pic}^0(C)\longrightarrow 0$ induces 
an embedding $\Phi_L\hookrightarrow U_L^{{\rm ab}}$. One 
identifies $\Phi_L$ with the $\widehat{{\Bbb Z}}$-module 
of the  $\widehat{{\Bbb Z}}(1)$-valued functions on $C(k)$ 
modulo the constant ones. 

The next step is to get a description of $\Phi_L$ in 
terms of the Galois groups. Clearly, $U_{k(x)}^{{\rm ab}}=\Phi_{k(x)}$. 
\begin{lemma} \begin{itemize} \item If $U$ is an open compact 
subgroup in $G$ then $N_G(U)=N_{G^{\circ}}(U)$. 

If, moreover, $N_G(U)/U$ has no abelian subgroups of 
finite index then $U=U_{k(x)}$ for some $x\in F-k$. 
\item For any $x\in L-k$ the transfer $U^{{\rm ab}}_{k(x)}
\longrightarrow U^{{\rm ab}}_L$ factors through $\Phi_L$. 
\item The span of images of the transfers 
$U^{{\rm ab}}_{k(x)}\longrightarrow U^{{\rm ab}}_L$ for all $x\in L-k$ 
is dense in $\Phi_L$. \end{itemize} \end{lemma}
{\it Proof.} \begin{itemize} \item By Proposition \ref{main-psha}, $U=U_L$ for 
a field $L$ finitely generated over $k$. Then the group $N_G(U_L)/U_L$ 
coincides with the group of automorphisms of the field $L$ over $k$. 
As the automorphism groups of projective curves of genus $>1$ are 
finite, if $L$ is isomorphic to the function 
field of such a curve, then the normalizer of $U$ in $G$ is compact. 
As the automorphism groups of elliptic curves are generated by 
elements of order $\le 4$ and contain abelian subgroups of index 
$\le 6$, if $L$ is isomorphic to the function field of such a curve, 
then the normalizer of $U$ in $G$ is generated by its compact 
subgroups. This implies that if $N_G(U)/U$ has no abelian subgroups 
of finite index then $L$ should be the function field of a rational 
curve. As the automorphism group of the rational curve is generated 
by involutions, the normalizer of $U$ in $G$ is generated by its 
compact subgroups. 
\item The transfer is induced by the norm homomorphism 
$L^{\times}/k^{\times}\stackrel{{\rm Nm}_{L/k(x)}}%
{-\!\!\!-\!\!\!-\!\!\!\longrightarrow}k(x)^{\times}/k^{\times}$, 
which is the restriction of the push-forward map 
${\rm Div}^0(C)\stackrel{x_{\ast}}{\longrightarrow}
{\rm Div}^0({\Bbb P}^1)$. 

Since $k(x)^{\times}/k^{\times}=
{\rm Div}^0({\Bbb P}^1)$, the transfer factors through $\Phi_L$. 
\item Each point $p$ of a smooth projective model $C$ 
of $L$ over $k$ is a difference of very ample effective divisors 
on $C$. These divisors themselves are zero-divisors of some 
rational functions, i.e., there are surjective morphisms 
$x,y:C\longrightarrow{\Bbb P}^1$ and a point $0\in{\Bbb P}^1$ such 
that $x^{-1}(0)-y^{-1}(0)=p$. Then 
$\delta_p=x^{\ast}\delta_0-y^{\ast}\delta_0:C(k)
\longrightarrow\widehat{{\Bbb Z}}(1)$ 
is a $\delta$-function of the point $p$ of $C$. As the span of 
$\delta$-functions is dense in the group 
$\Phi_L$, we are done. \hfill $\Box$ \end{itemize} 
\vspace{5mm}

For a point of $C(k)$ its {\sl decomposition subgroup} in 
$\Phi_L\subset U_L$ consists of all functions supported on it. 
In the case $L=k(x)$ the decomposition subgroups in $U_{k(x)}^{{\rm ab}}$ 
are parametrized by the set (which is isomorphic to ${\Bbb P}^1(k)$) 
of parabolic subgroups $P$ in $N_GU_{k(x)}/U_{k(x)}$. 
The subgroup $D_P$ consists of elements in $U_{k(x)}^{{\rm ab}}$ 
fixed under the adjoint action of $P$. 
Clearly, $D_P\cong\widehat{{\Bbb Z}}(1)$. 

Each inclusion of subgroups $U_L\subset U_{k(x)}$ induces 
a homomorphism $U_L^{{\rm ab}}\longrightarrow U_{k(x)}^{{\rm ab}}$. Consider 
the evident homomorphism $U_L^{{\rm ab}}\stackrel{\varphi_L}{\longrightarrow}
\prod_{x\in L-k}U_{k(x)}^{{\rm ab}}$. For any non-zero element $h$ of 
the group $U_L^{{\rm ab}}$, considered as a homomorphism from the group 
$L^{\times}$, there is an element $x\in L^{\times}$ with $h(x)\neq 0$, so 
the image of $h$ in $U_{k(x)}^{{\rm ab}}$ is non-zero, and thus, 
$\varphi_L$ is injective. 

To construct decomposition subgroups for an arbitrary $L$, consider 
such a subgroup $D\cong\widehat{{\Bbb Z}}$ in the target of $\varphi_L$ 
that its projection to each of $U_{k(x)}^{{\rm ab}}$ is of finite index 
in some decomposition subgroup. Then our nearest goal is to show 
that the set of decomposition subgroups in $U_L^{{\rm ab}}$ coincides 
with the set of maximal subgroups among $\Phi_L\cap\varphi^{-1}_L(D)$. 
\begin{lemma}[ = Lemma 5.2 of \cite{ab-subgroups} = 
Lemma 3.4$\vphantom{4}'$ of \cite{bo-two}]  \label{Lemma 5.2} 
Suppose that $f$ is such a function on a projective space ${\Bbb P}$ over 
an infinite field that the restriction of $f$ to each projective line in 
${\Bbb P}$ is constant on the complement to a point on it. 

Then $f$ is a flag function, i.e., there is a filtration 
$P_0\subset P_1\subset P_2\subset\dots$ of ${\Bbb P}$ 
by projective subspaces such that $f$ is constant on $P_0$ and on all 
strata $P_{j+1}-P_j$. \hfill $\Box$ \end{lemma}

\begin{lemma} \label{vanishing} For any smooth projective 
curve $C$ there is a constant $N=N(C)$ such that for any 
${\cal L},{\cal L}'\in{\rm Pic}^{\ge N}(C)$ the natural map 
$\Gamma(C,{\cal L})\otimes\Gamma(C,{\cal L}')\longrightarrow
\Gamma(C,{\cal L}\otimes{\cal L}')$ is surjective. \end{lemma}
{\it Proof.} Fix an invertible sheaf ${\cal L}_0$ on $C$ of 
degree 1. By Serre vanishing theorem, there is such an integer $N'$ 
that the sheaf 
$\left({\cal L}_0\boxtimes{\cal L}_0\right)^{\otimes N'}(-\Delta)$ 
on $C\times C$ is generated by its global sections, and therefore, 
for any ${\cal L},{\cal L}'\in{\rm Pic}^{>N'}(C)$ the sheaf 
$({\cal L}\boxtimes{\cal L}')(-\Delta)$ is ample. 
Let $N=N'+2g$. Then by Kodaira vanishing theorem, for any 
${\cal L},{\cal L}'\in{\rm Pic}^{\ge N}(C)$ the short exact 
sequence $$0\longrightarrow({\cal L}\boxtimes{\cal L}')(-\Delta)
\longrightarrow{\cal L}\boxtimes{\cal L}'\longrightarrow
{\cal L}\otimes{\cal L}'\longrightarrow 0$$ 
of sheaves on $C\times C$ induces a surjection 
$\Gamma(C,{\cal L})\otimes_k\Gamma(C,{\cal L}')\longrightarrow
\Gamma(C,{\cal L}\otimes{\cal L}')$. \hfill $\Box$

\begin{lemma} If $\varphi^{-1}_L(D)$ is in $\Phi_L$ then it is a 
subgroup in a decomposition subgroup in $U_L^{{\rm ab}}$.\end{lemma}
{\it Proof.} Let $f\in\varphi^{-1}_L(D)\cap\Phi_L$, i.e., 
$f:C(k)\longrightarrow\widehat{{\Bbb Z}}(1)$ for a smooth projective 
model $C$ of $L$ over $k$, and for any very ample invertible sheaf 
${\cal L}$ on $C$ restrictions of the induced function 
$f:|{\cal L}|\longrightarrow\widehat{{\Bbb Z}}(1)$ 
to projective lines in $|{\cal L}|$ are ``$\delta$-functions'' on them. 
Then, by Lemma \ref{Lemma 5.2}, $f$ is a flag function. Therefore, the 
function $\widehat{f}:|{\cal L}|^{\vee}\longrightarrow
\widehat{{\Bbb Z}}(1)$ given by 
$H\longmapsto f(\mbox{general point of $H$})$ is a ``$\delta$-function''. 

Let $g$ be the genus of $C$. Consider the composition 
$\widehat{f}_{{\cal L}}:C(k)\longrightarrow|{\cal L}|^{\vee}
\stackrel{\widehat{f}}{\longrightarrow}\widehat{{\Bbb Z}}(1)$. 
It takes $x$ to $$f(x)+f(\mbox{general point of $|{\cal L}(-x)|$}).$$ 
Since it is a ``$\delta$-function'', and all the hyperplanes 
$x+|{\cal L}(-x)|$ in 
${\cal L}$ are pairwise distinct, there are such functions 
$b_0:{\rm Pic}^{> 2g}(C)\longrightarrow\widehat{{\Bbb Z}}(1)$ and 
$a:{\rm Pic}^{> 2g}(C)\longrightarrow C(k)$ that 
$$f(x)+f(\mbox{general point of $|{\cal L}(-x)|$})=b_0({\cal L})
\delta_{x,a({\cal L})}+b_1({\cal L}),$$ 
where $b_1:{\rm Pic}^{\ge 2g}(C)\longrightarrow\widehat{{\Bbb Z}}(1)$ 
is the function 
${\cal L}\longmapsto f(\mbox{general point of $|{\cal L}|$})$. Then 
$$f(x)=b_0({\cal L})\delta_{x,a({\cal L})}+b_1({\cal L})-b_1({\cal L}(-x)).$$

By Lemma \ref{vanishing}, for any ${\cal L},{\cal L}'
\in{\rm Pic}^{\ge N}(C)$ the image of the map 
$|{\cal L}|\times|{\cal L}'|\longrightarrow|{\cal L}\otimes{\cal L}'|$ 
of summation of divisors is not contained in any hyperplane in 
$|{\cal L}\otimes{\cal L}'|$. Then a sum of a general divisor in 
$|{\cal L}|$ and a general divisor in $|{\cal L}'|$ is a general 
divisor in the linear system $|{\cal L}\otimes{\cal L}'|$, so one 
has $$b_1({\cal L}\otimes{\cal L}')=b_1({\cal L})+b_1({\cal L}'),$$ 
and therefore, for any sheaf ${\cal L}_0$ of degree zero one has 
$$b_1({\cal L}')+b_1({\cal L}_0\otimes{\cal L})=
b_1({\cal L})+b_1({\cal L}_0\otimes{\cal L}'),$$
so $b_2({\cal L}_0):=b_1({\cal L}_0\otimes{\cal L})-b_1({\cal L}):
{\rm Pic}^0(C)\longrightarrow\widehat{{\Bbb Z}}(1)$ does not depend on 
${\cal L}$. It is easy to see that $b_2$ is a homomorphism, which therefore 
should be zero, since ${\rm Pic}^0(C)$ is a divisible group. From this we 
conclude that $b_1({\cal L})=b_1(\deg{\cal L})$, and finally, 
$f(x)=b_0({\cal L})\delta_{x,a({\cal L})}+b_3({\cal L})$ is a 
$\delta$-function on $C(k)$, i.e., corresponds to a point of $C$, 
or to a decomposition subgroup in $U_L^{{\rm ab}}$. \hfill $\Box$

\section{Automorphisms of subgroups between $G^{\circ}$ and $G$}
\begin{lemma} \label{auto-f} \begin{enumerate} \item \label{auto-ff} 
Suppose that for a subgroup $H$ in $G$ containing $G^{\circ}$ 
(the restriction to $G^{\circ}$ of) a homomorphism 
$\lambda:H\longrightarrow G$ 
induces the identity map of the set ${\frak F}$ of compact open 
subgroups in $G$. Then $\lambda=id$. 
\item \label{auto-ffff} The centralizer of $G^{\circ}$ in 
$G_{F/{\Bbb Q}}$ is trivial. \end{enumerate} \end{lemma}
{\it Proof.} For any $\sigma\in H$ and any open compact subgroup 
$U$ one has $\sigma U\sigma^{-1}=\lambda(\sigma U\sigma^{-1})
=\lambda(\sigma)\lambda(U)\lambda(\sigma)^{-1}=
\lambda(\sigma)U\lambda(\sigma)^{-1}$, so $\sigma^{-1}\lambda(\sigma)$ 
belongs to the normalizer of each $U$. 

For a variety $X$ of dimension $n$ over $k$ without birational 
automorphisms and any $x\in F-k$ there is a subfield $L_x\subset F$ 
containing $x$ isomorphic to the function field of $X$. Then the 
normalizer of $U_{L_x}$ coincides with $U_{L_x}$, and the intersection 
of all $U_{L_x}$ is $\{1\}$, so $\sigma^{-1}\lambda(\sigma)=1$. On the 
other hand, if $\tau\in G_{F/{\Bbb Q}}$ normalizes $U_{k(x,P(x)^{1/2})}$ 
for all polynomials $P$ over $k$, then $\tau\in G_{F/k}$, and 
therefore, $\tau=1$. \hfill $\Box$ 

\vspace{4mm}
Let ${\frak F}$ be the set of compact open subgroups 
in $G^{\circ}$, and let ${\Bbb Q}(\chi)$ be the quotient of the 
free abelian group generated by ${\frak F}$ by the relations 
$[U]=[U:U']\cdot[U']$ for all $U'\subset U$. As the intersection 
of a pair of a compact open subgroups in $G$ is a subgroup of finite 
index in both of them, ${\Bbb Q}(\chi)$ is a one-dimensional 
${\Bbb Q}$-vector space. The group $G$ acts on it by the 
conjugations. Let $\chi$ be the character of this representation of 
$G$. 

One can get an explicit formula for $\chi$ as follows. Fix a subfield 
$L$ of $F$ finitely generated and of transcendence degree $n$ over $k$. 
Then for any $\sigma\in G$ one has 
$[U_L]=[L\sigma(L):L]\cdot[U_{L\sigma(L)}]$ and 
$[U_{\sigma(L)}]=[L\sigma(L):\sigma(L)]\cdot[U_{L\sigma(L)}]$, 
and therefore, 
$\chi(\sigma)=\frac{[L\sigma(L):\sigma(L)]}{[L\sigma(L):L]}$. 
This implies that $\chi:G\longrightarrow{\Bbb Q}^{\times}_+$ is 
surjective, and its restriction to $G^{\circ}$ is trivial.

\begin{theorem} \label{auto-gen} Let $n=1$, $H$ be a subgroup in 
$G$ containing $G^{\circ}$, and $N_{G_{F/{\Bbb Q}}}(H)$ its 
normalizer in $G_{F/{\Bbb Q}}$. Then $N_{G_{F/{\Bbb Q}}}(H)
\subseteq N_{G_{F/{\Bbb Q}}}(G)=\{\mbox{{\rm automorphisms of} 
$F$ {\rm preserving} $k$}\}$, and the adjoint action of 
$N_{G_{F/{\Bbb Q}}}(H)$ on $H$ induces an isomorphism from 
$N_{G_{F/{\Bbb Q}}}(H)$ to the group of continuous open 
automorphisms of $H$. If $H\supseteq\ker\chi$ then 
$N_{G_{F/{\Bbb Q}}}(H)=N_{G_{F/{\Bbb Q}}}(G^{\circ})$. \end{theorem} 
{\it Proof.} For each $U\in{\frak F}$ let ${\rm Div}_U^+$ be the 
free abelian semi-group, whose generators are decomposition 
subgroups in $U^{{\rm ab}}$, and for each integer $d\ge 2$ let 
$${\frak Gr}^{(d)}_U=\{U_L\supset 
U\quad |\quad [U_L:U]=d,L\cong k(t)\}\subset{\frak F}.$$ 
For a smooth projective model $C$ of $F^U$ the set ${\frak Gr}^{(d)}_U$ 
is in bijection with the disjoint union of Zariski-open subsets in 
Grassmannians $$\coprod\limits_{{\cal L}\in{\rm Pic}^d(C)}\left(
Gr(1,|{\cal L}|)-\bigcup_{x\in C(k)}Gr(1,x+|{\cal L}(-x)|)\right).$$ 

One can define~: \begin{itemize} 
\item an {\sl ``invertible sheaf of degree $d$ without base points''} 
${\cal L}$, as a subset of ${\frak Gr}^{(d)}_U\subset{\frak F}$ 
consisting of elements equivalent under the relation generated by 
$U_1\sim_UU_2$ if there are decomposition subgroups 
$D_a\subset U^{{\rm ab}}_1$ and $D_b\subset U^{{\rm ab}}_2$ 
such that their preimages in $U^{{\rm ab}}$ contain the same 
collections of decomposition subgroups with the same indices of their
images in $D_a$ and $D_b$; 
\item the {\sl ``linear system''} $|{\cal L}|$, as the set of maximal 
collections of elements of ${\cal L}$ ``intersecting at a single point'', 
i.e., as the subset of the free abelian semi-group ${\rm Div}^+_U$; 
\item a {\sl ``line presented in ${\cal L}$''} in $|{\cal L}|$, 
as an element of 
${\cal L}\subset{\frak Gr}^{(d)}_U$, considered as a subset in $|{\cal L}|$; 
\item an arbitrary {\sl ``line''} in $|{\cal L}|$, as a subset in 
$|{\cal L}|$ of type $D+l$, where $D\in{\rm Div}^+_U$ and $l$ is a 
line presented in the sheaf ${\cal L}(-D)$ without base points; 
\item an ``$s$-{\sl subspace}'' in $|{\cal L}|$, as the union of 
all lines passing through a given point in $|{\cal L}|$ and intersecting 
a given ``$(s-1)$-subspace'' in $|{\cal L}|$. \end{itemize} 

Now we remark that for any sufficiently big $d$ and any sheaf 
${\cal L}\subset{\frak Gr}^{(d)}_U$ the set $C_U$ of decomposition 
subgroups in $U^{{\rm ab}}$ can be {\sl canonically} identified with 
the subset of $|{\cal L}|^{\vee}$ consisting of those hyperplanes in 
$|{\cal L}|$ that each line on each of them is ``absent in ${\cal L}$''. 
As $|{\cal L}|^{\vee}$ has a canonical structure of a projective 
space (but not of a projective space over $k$), this gives 
us a {\sl canonical} structure of a scheme on $C_U$. Let $\kappa_U$ be 
the function field of $C_U$. 

Clearly, $\lambda(G^{\circ})=G^{\circ}$ and the restriction of 
$\lambda$ to $G^{\circ}$ induces a bijection 
${\frak Gr}^{(d)}_U\stackrel{\sim}{\longrightarrow}
{\frak Gr}^{(d)}_{\lambda(U)}$ for each $d\ge 2$, and for 
any sheaf ${\cal L}\subset{\frak Gr}^{(d)}_U$ it induces a map 
$|{\cal L}|\longrightarrow|\lambda({\cal L})|$ which transforms 
subspaces into subspaces (of the same dimension), i.e., a collineation. 
As $\lambda$ induces a collineation 
$|{\cal L}|^{\vee}\stackrel{\sim}{\longrightarrow}|\lambda({\cal L})|^{\vee}$, 
the fundamental theorem of projective geometry (see, e.g., 
\cite{artin}) implies that such $\lambda$ induces an isomorphism 
$C_U\stackrel{\sim}{\longrightarrow}C_{\lambda(U)}$ of schemes over 
${\Bbb Q}$. 
This isomorphism does not depend on $d$ and ${\cal L}$, since it 
determines the collineations $|{\cal L}'|\stackrel{\sim}{\longrightarrow}
|\lambda({\cal L}')|$ for all ${\cal L}'\subset{\frak Gr}^{(d')}_U$. 
Denote by $\sigma_U$ the induced isomorphism 
$\kappa_{\lambda(U)}\stackrel{\sim}{\longrightarrow}\kappa_U$. 

For each subgroup $U'$ of finite index in $U$ the natural map 
$C_{U'}\longrightarrow C_U$ is a morphism of schemes, and in particular, 
$\kappa_U$ is naturally embedded into $\kappa_{U'}$. The group 
$G^{\circ}$ acts on the field $\lim\limits_{_U\longrightarrow}\kappa_U$. 
By Lemma \ref{auto-f} (\ref{auto-ffff}), the centralizer of $G^{\circ}$ 
in $G$ is trivial, and therefore, there is a unique isomorphism 
$\lim\limits_{_U\longrightarrow}
\kappa_U\stackrel{\sim}{\longrightarrow}F$ commuting with 
the $G^{\circ}$-action. Since the diagram 
$$\begin{array}{ccc} C_{U'} & \longrightarrow & C_{\lambda(U')} \\
\downarrow & & \downarrow \\
C_U & \longrightarrow & C_{\lambda(U)} \end{array}$$ 
commutes, the restriction of $\sigma_{U'}$ to $\kappa_U$ coincides 
with $\sigma_U$, and finally, we get an automorphism $\sigma$ of $F$ 
induced by $\lambda$. As $k$ is the only maximal algebraically closed 
subfield in its arbitrary finitely generated extension, $\sigma$ 
induces an automorphism of $k$, and therefore, normalizes $G^{\circ}$. 

Then the restriction to $G^{\circ}$ of ${\rm ad}(\sigma)\circ\lambda$ 
acts trivially on all of ${\frak Gr}^{(d)}_{U'}$. As any open compact 
subgroup is an intersection of elements of ${\frak Gr}^{(d)}_{U'}$ for 
$d$ big enough and $U'$ small enough, ${\rm ad}(\sigma)\circ\lambda$ 
acts on ${\frak F}$ also trivially. By Lemma \ref{auto-f} (\ref{auto-ff}), 
this implies that $\lambda={\rm ad}(\sigma^{-1})$. \hfill $\Box$ 

\vspace{4mm}

{\sc Remark.} If $k$ is countable then the inverse of 
any continuous automorphism as in the statement of 
Theorem \ref{auto-gen} is automatically continuous~:
\begin{lemma} \label{aut-open} If $k$ is countable, and 
$U\stackrel{\lambda}{\longrightarrow}U'$ is a continuous surjective 
homomorphism 
of open subgroups in $G_{F/k}$ and $G_{F'/k'}$ then the image in $U'$ 
of an open subset in $U$ is open. \end{lemma} 
{\it Proof.} Let $U_L\subset U$ be an open compact subgroup. Then 
$U/U_L$ is a countable set surjecting onto the set $U'/\lambda(U_L)$. 
By Proposition \ref{main-psha}, for the subfield $L'=F^{\lambda(U_L)}$ 
one has $\overline{L'}=F$. If $\lambda(U_L)$ is not open then $L'$ is 
not finitely generated over $k'$, and therefore, $U'/\lambda(U_L)$ 
is not countable. \hfill $\Box$

\end{document}